\documentclass[leqno,12pt]{article}
\usepackage{latexsym}
\usepackage{amsmath}
\usepackage{amssymb}
\usepackage{bbm}
\usepackage{color}
\usepackage{pst-all}
\usepackage[authoryear]{natbib}
\usepackage{graphicx}
\usepackage{epsfig}

 0

\newtheorem{thm}{Theorem}[section]
\newtheorem{prop}[thm]{Proposition}

\newtheorem{lem}[thm]{Lemma}

\makeatletter\@addtoreset{equation}{section}\makeatother

\newcommand{\ba}{\begin{array}}
\newcommand{\ea}{\end{array}}
\newcommand{\beqohne}{\begin{eqnarray*}}
\newcommand{\eeqohne}{\end{eqnarray*}}
\newcommand{\beohne}{\begin{equation*}}
\newcommand{\eeohne}{\end{equation*}}

\def\proof{\noindent{\bf Proof:}\hskip10pt}
\def\QED{\hfill $ \Box $\\}

\begin{document}

\title{Moderate deviations for spectral measures of random matrix ensembles}

\author{Jan Nagel}

\maketitle
\begin{abstract}
In this paper we consider the (weighted) spectral measure $\mu_n$ of a $n\times n$ random matrix, distributed according to a classical Gaussian, Laguerre or Jacobi ensemble, and show a moderate deviation principle for the standardised signed measure $\sqrt{n/a_n}(\mu_n -\sigma)$. The centering measure $\sigma$ is the weak limit of the empirical eigenvalue distribution and the rate function is given in terms of the $L^2$-norm of the density with respect to $\sigma$. The proof involves the tridiagonal representations of the ensembles which provide us with a sequence of independent random variables and a link to orthogonal polynomials.
\end{abstract}

\medskip

\section{Introduction}

In the theory of random matrices the three central distributions on the set of self-adjoint matrices are the Gaussian, the Laguerre and the Jacobi ensemble. For these distributions, a lot is known about the asymptotic properties of the empirical eigenvalue distribution
\begin{align*} 
\hat\mu_n = \frac{1}{n} \sum_{i=1}^n  \delta_{\lambda_i} ,
\end{align*}
when $\lambda_1,\dots ,\lambda_n$ are the eigenvalues of the matrix under consideration and $\delta_x$ denotes the Dirac measure in $x$. The random measure $\hat\mu_n$ converges almost surely weakly to a limit distribution $\sigma$ with compact support. For the Gaussian ensemble, this is the famous semicircle law (\cite{wigner1958}). In the Laguerre case and Jacobi case, the limit distributions are the Marchenko-Pastur law and the semicircle law (see \cite{marpas1967} and \cite{col2005}). An overview and more recent results can be found in the monographs of \cite{mehta2004}, \cite{forrester2010} and \cite{agz2010}.

Recently, a number of paper deal with the weighted version of the empirical eigenvalue distribution, the spectral measure defined as
\begin{align} \label{wspectralmeasure}
\mu_n = \sum_{i=1}^n w_i \delta_{\lambda_i} ,
\end{align}
with weights $w_i = |\langle u_i,e_1\rangle |^2$ given by the squared modulus of the top entries of the eigenvectors $u_1,\dots ,u_n$, see the paper of \cite{gamloz2004}, \cite{gamrou2009b,gamrou2009a}.
This weighted measure encodes some information about the eigenvectors and appears in the general spectral theorem. 
For the general formulation suppose $A$ is a linear, continuous and self-adjoint operator on a Hilbert space $\mathcal{H}$ with a cyclic vector $h \in \mathcal{H}$ (i.e., the span of $h, Ah, A^2h\dots $ is $\mathcal{H}$). The spectral theorem yields the existence of a spectral measure $\mu$ on some measureable space $(X, \mathcal{X})$, a unitary isometry $U: L^2(\mu) \rightarrow \mathcal{H}$ and a bounded measureable function $\phi$ on $X$ such that
\begin{align*}
\left( U^* A U f \right) (x) = \phi (x) f(x)
\end{align*}
for each $f \in L^2(\mu)$ (see \cite{dunfschw1963}). This statement is sometimes abbreviated to ``every Hermitian operator is unitarily equivalent to a multiplication''. The isometry $U$ maps a polynomial $p(x)$ to $p(A) h$, which implies that the moments of $\mu$ satisfy
\begin{align*} 
\int x^k \mu(dx) = \langle 1,x^k \rangle_{L^2} = \langle h,A^k h \rangle 
\end{align*}
for all $k\geq 1$. If $\mathcal{H}$ is $n$-dimensional and $h=e_1$, then this gives the moments of \eqref{wspectralmeasure}. By the functional calculus, $\int fd\mu_n = f(A)_{1,1}$ for suitable test functions $f$, so the behaviour of $\mu_n$ is related to matrix elements of functions of random matrices, which were studied by \cite{lytpast2009,lytpast2009b} and \cite{pizzo2012}.

\smallskip

We are concerned with moderate deviation properties of the spectral measure. To be self-contained, we recall the definition of large and moderate deviations.
Let $U$ be a topological Hausdorff space, let $\mathcal{I}: U \rightarrow [0,\infty]$ be a lower semicontinuous function and $a_n$ a sequence of positive real numbers with $a_n\to \infty$. We say that a sequence $(P_{n})_n$ of probability measures on $(U,\mathcal{B}(U))$ satisfies a large deviation principle (LDP) with speed $a_n$ and
rate function $\mathcal{I}$  if:
\begin{itemize}
\item[(i)] For all closed sets $F \subset U$:
\begin{align*}
\limsup_{n\rightarrow\infty} \frac{1}{a_n} \log P_{n}(F)\leq -\inf_{x\in F}\mathcal{I}(x)
\end{align*}
\item[(ii)] For all open sets $O \subset U$:
\begin{align*}
\liminf_{n\rightarrow\infty} \frac{1}{a_n} \log P_{n}(O)\geq -\inf_{x\in O}\mathcal{I}(x)
\end{align*}
\end{itemize}
The rate function $\mathcal{I}$ is good if its level sets
$\{x\in U |\ \mathcal{I}(x)\leq a\}$ are compact for all $a\geq 0$.
We say that a sequence of $U$-valued random variables satisfies an LDP if their distributions satisfy an LDP. 
If $a_n \to \infty$, but $a_n/n \to 0$, we say that  $(P_{n})_n$ satisfies a moderate deviation principle (MDP). Moderate deviations are concerned with the convergence of measure at a scaling between the law of large numbers and the central limit theorem. On this scale, the rate function is identical for all possible speeds.

The empirical eigenvalue distribution $\hat\mu_n$ satisfies a large deviation principle with speed $n^2$ and rate function related to Voiculescu's entropy, see \cite{bengui1997}, \cite{hiaipetz1998} and \cite{hiaipetz2006}. The spectral measure $\mu_n$ has the same weak limit $\sigma$ as the empirical distribution, the large deviation properties however, are different: the random weights reduce the speed to $n$. In the case of the Gaussian ensemble, the rate function is explicit and involves the Kullback-Leibler distance and a sum over atoms outside the limit support (\cite{gamrou2009b}). In several cases, moderate deviations of empirical eigenvalue distributions were proven. \cite{demguizei2003} considered the perturbation of a Gaussian matrix and \cite{doereich2010,doereich2011} showed moderate deviations for Wigner matrices generalising the Gaussian ensemble. 
Our main result is Theorem \ref{haupt}, a moderate deviation principle for the rescaled measure $\sqrt{n/a_n}(\mu_n -\sigma)$, when $\mu_n$ is the random spectral measure of one of the three classical ensembles.

The proof of Theorem \ref{haupt} consists of two main steps. The first one is a moderate deviation principle for the moments of the rescaled measure. The key part is a transformation to new random variables which are independent and related to the polynomials orthogonal with respect to the spectral measure. In the second step we apply the contraction principle to lift the result into the space of measures. Due to this method, the MDP will be in the topology induced by the convergence of moments.

This paper is structures as follows: in Section 2 we introduce the classical ensembles and corresponding spectral measures and we state our results. The transformation in the context of orthogonal polynomials is formulated in Section 3. Finally, Section 4 contains the proofs.

\medskip

\section{Random matrix ensembles, results}

The central distribution of Hermitian matrices and a generalisation of the scalar normal distribution is the Gaussian ensemble. Following the threefold way of \cite{dyson1962b} it comes in a real, complex and quaternion version and is motivated by several physical observations (see \cite{mehta2004}). In the real case the Gaussian orthogonal ensemble (GOE) is the distribution of a $n\times n$ matrix with centered Gaussian entries such that the random variables on and above the diagonal are independent, the diagonal entries have variance $\frac{2}{n}$ and the off-diagonal entries have variance $\frac{1}{n}$. For the complex Gaussian unitary ensemble (GUE), the off diagonal entries are complex Gaussian distributed with variance $\frac{1}{2n}$, the real diagonal entries have variance $\frac{1}{n}$. The quaternion matrices can be constructed in a similar way. In all three cases, the density of the eigenvalues can be written as
\begin{align} \label{ewge}
f_G(\lambda) = c_G | \Delta (\lambda )|^{\beta} \prod_{i=1}^n e^{-\frac{\beta n}{4}\lambda_i^2},
\end{align}
with $| \Delta (\lambda )|=\prod_{i< j} |\lambda_i-\lambda_j|$ the Vandermonde determinant, $\beta$ the real dimension of the entries and $c_G$ a normalisation constant.

The Laguerre ensemble (or Wishart ensemble) appears as the distribution of covariance estimators in multivariate statistics (\cite{wishart1928}). It can be constructed as the square $XX^*$ of a real, complex or quaternion Gaussian matrix $X$, the eigenvalue density is given by 
\begin{align} \label{ewle}
f_L(\lambda) = c_L^\gamma | \Delta (\lambda )|^{\beta} \prod_{i=1}^n \lambda_i^{\gamma} e^{-\frac{\beta n}{2} \lambda_i} \mathbbm{1}_{\{ \lambda_i>0\} },
\end{align}
where $\gamma>-1 $ is a parameter determined by the size of $X$.

The third distribution is the Jacobi ensemble with applications in multivariate analysis of variance (MANOVA, see \cite{muirhead1982}) and eigenvalue density
\begin{align} \label{ewje}
f_J(\lambda) = c_J^{\gamma,\delta} | \Delta (\lambda )|^{\beta} \prod_{i=1}^n \lambda_i^{\gamma} (1-\lambda_i)^{\delta} \mathbbm{1}_{\{ 0<\lambda_i<1 \} } .
\end{align}
and with parameters $\gamma ,\delta >-1$. The corresponding matrix model can be written as the ratio of two Laguerre matrices.

It is a common feature of the three ensembles that the distribution is invariant under orthogonal (or unitary or symplectic) conjugations. As a consequence, the matrix of eigenvalues is Haar distributed on the corresponding compact group and independent of the eigenvalues, see \cite{dawid1977}. Consequently, the first unit vector is almost surely cyclic and the weight vector containing the square moduli of the top entries of the eigenvectors
\begin{align*} 
w = (w_1,\dots ,w_n) = \left( |\langle u_1,e_1\rangle |^2,\dots ,|\langle u_n,e_1\rangle |^2 \right) 
\end{align*}
follows a Dirichlet distribution Dir($\frac{\beta}{2}$) with density
\begin{align} \label{dirdensity}
\frac{\Gamma (\tfrac{n\beta}{2} )}{\Gamma (\tfrac{\beta}{2})^n} w_1^{\beta/2-1} \dots w_n^{\beta/2-1} \mathbbm{1}_{\{ w_1+\dots + w_n=1, w_i> 0 \} } 
\end{align}
(\cite{johnson1976}). 
The distribution of the eigenvalues and the Dirichlet distribution of the weights $w$ are well-defined not only for $\beta \in \{1,2,4\}$ but for all values $\beta >0$ and we can consider random spectral measures of the three ensembles in the general case. To be precise, we say that 
\begin{align*}
\mu_n = \sum_{i=1}^n w_i \delta_{\lambda_i}
\end{align*}
is a spectral measure of a classical ensemble if the distribution of the support points is given by \eqref{ewge}, \eqref{ewle} or \eqref{ewje} independent of the weights, which have density \eqref{dirdensity}.

\medskip

The eigenvalue distributions defined here are normalised such that the empirical eigenvalue distribution $\hat\mu_n$ has a weak limit almost surely. Moreover, they satisfy a large deviation principle with speed $n^2$ and good rate functions related to Voiculescu's entropy. For the Gaussian ensemble, the equilibrium measure and thus the limit of $\hat\mu_n$ is Wigner's famous semicircle law
\begin{align}\label{semicircle}
\rho(dx) = \frac{1}{2\pi} \sqrt{4-x^2} dx .
\end{align}
In the Laguerre case, the limit is the Marchenko-Pastur distribution
\begin{align} \label{marchpast}
\eta(dx) = \frac{\sqrt{x(4-x)}}{2\pi x} \mathbbm{1}_{ \{ 0<x<4 \} } dx
\end{align}
and the eigenvalue distribution of the Jacobi ensemble converges to the arcsine law
\begin{align} \label{arcsine}
\nu (dx) = \frac{1}{\pi \sqrt{x(1-x)}} \mathbbm{1}_{ \{ 0<x<1 \} } dx .
\end{align}
The same convergence holds for the weighted versions, the spectral measures $\mu_n$. More precisely, suppose that the expected moments of $\hat\mu_n$ converge, then bounding the moments of the Dirichlet weights, we can show
\begin{align*}
\left| m_k(\mu_n) - m_k (\hat{\mu}_n)\right| = O_P\left( \frac{1}{\sqrt{n}}  \right)
\end{align*}
 for all $k$ as $n\to \infty$, where 
\begin{align*}
m_k(\mu) = \int x^k \mu(dx)
\end{align*}
denotes the $k$-th moment of a measure $\mu$. Interestingly, proving convergence of expected moments is precisely the way a number of limit theorems in the framework of free probability are proven, where a comprehensive overview is given by \cite{nicspe2006}.

\medskip

In order to present the results in a unified way and to avoid repetitions, we denote by $\mu_n$ the spectral measure of one of the three ensembles and the limit law by $\sigma$. We regard the standardised measure $\mu_n-\sigma$ as a random element of the set \begin{align*}
\mathcal{M}_0 = \left\{ \mu | \mu \mbox{ finite signed measure on } \mathbb{R} \mbox{ with compact support, } \mu(\mathbb{R})=0 \right\} .
\end{align*}
We endow the space $\mathcal{M}_0$ with the topology of convergence of moments, that is, $\mu_n \to \mu$ if for all moments $m_k(\mu_n) \to m_k(\mu)$ holds. The family of monomials $x^k$ is separating on $\mathcal{M}_0$, in the sense that for $\mu \neq \nu$ there is a $k$ such that $\int x^kd\mu \neq \int x^k d\nu$. A metric for this topology is then given by
\begin{align} \label{metric}
d_\mathbf{m} (\mu,\nu) =\sum_{k=0}^\infty 2^{-k} \frac{\left| \int x^k d(\mu-\nu) \right|}{1+\left| \int x^k d(\mu-\nu) \right|} .
\end{align} 
Note that on $\mathcal{M}_0$, the topology of convergence of moments is neither stronger nor weaker than the topology of weak convergence. The $\sigma$-algebra is then the usual Borel algebra. The following theorem is our main result.

\medskip

\begin{thm}\label{haupt}
Let $\mu_n$ be a random spectral measure of a classical ensemble with weak limit $\sigma$. Then $(\sqrt{n/a_n}(\mu_n - \sigma))_n$ satisfies the moderate deviation principle in $\mathcal{M}_0$ with speed $a_n$ and good rate function
\begin{align*}
\mathcal{I}(\mu) = \frac{\beta}{4} \int \left( \frac{\partial \mu}{\partial \sigma }\right)^2 d\sigma ,
\end{align*}  
where $\frac{\partial \mu}{\partial \sigma }$ is the signed density with respect to $\sigma$. If $\mu$ is not absolutely continuous with respect to $\sigma$, we have $\mathcal{I}(\mu)=\infty$.
\end{thm}

\medskip

%

For the first part of the proof, we consider the behaviour of moments of $\mu_n$. 
Motivated from the theory of random moment sequences, \cite{detnag2012} obtained central limit theorems for the moment vectors of the classical ensembles. Let
$m^{(k)}(\mu) = (m_1(\mu) ,\dots ,m_k(\mu))^T$ denote the vector of the first $k$ moments. Then the convergence in distribution
\begin{align} \label{clt}
\sqrt{\frac{\beta n}{2}}  (m^{(k)}(\mu_n) - m^{(k)}(\sigma)) \xrightarrow[n \rightarrow \infty ]{\mathcal{D}} \mathcal{N}_k(0,D_kD_k^T),
\end{align}
holds, where $D_k$ is a lower triangular matrix depending on the ensemble under consideration. For the Gaussian ensemble, 
\begin{align}\label{covgau}
(D_k^{(G)})_{i,j} = \binom{i}{\frac{i-j}{2}} - \binom{i}{\frac{i-j}{2} - 1}.
\end{align}
for $i\geq j$ and $i+j$ even and $(D_k^{(G)})_{i,j}=0$ otherwise. As usual we let $\binom{i}{-1}=0$. The asymptotic covariance for the Laguerre ensemble satisfies
\begin{align}\label{covlag}
(D_k^{(L)})_{i,j} = \binom{2i}{i-j} - \binom{2i}{i-j-1} \qquad j \le i 
\end{align}
and in the Jacobi case we have
\begin{align} \label{covjac}
(D_k^{(J)})_{i,j} = \frac{2^{-2i+1}}{\sqrt{2}} \binom{2i}{i-j} \qquad j \le i .
\end{align}
From the central limit theorem, we are able to prove  a moderate deviation principle on the level of moments. 

\medskip

\begin{thm}\label{step1}
The sequence of moment vectors $\sqrt{n/a_n}(m^{(k)}(\mu_n) - m^{(k)}(\sigma))$ satisfies a moderate deviation principle with speed $a_n$ and good rate function
\begin{align*}
\mathcal{I}_m(m)= \frac{\beta}{4} ||D_k^{-1} m||^2 ,
\end{align*}
where $D_k$ is the matrix in \eqref{clt}.
\end{thm}

\medskip

In order to apply the contraction principle, we need an interpretation of $D_k^{-1}$ in terms of $\sigma$. As it turns out, the covariance matrix in \eqref{clt} has a very interesting structure.

\medskip

\begin{lem}\label{step2}
Let $\mu_n$ be the spectral measure of a classical ensemble with weak limit $\sigma$, then the covariance matrix in \eqref{clt} has the representation
\begin{align*}
(D_kD_k^T)_{i,j} = m_{i+j}(\sigma) - m_i(\sigma) m_j(\sigma) .
\end{align*}
\end{lem}

\medskip

\section{Orthogonal polynomials, tridiagonal representations}

Let $\mu$ be a measure on the real line with compact support, and let $p_0(x),p_1(x),\dots$ be the monic polynomials orthogonal with respect to $\mu$, i.e.,
\begin{align*}
\int p_k(x)p_l(x) \mu(dx) = c_k \delta_{k,l}
\end{align*}
and the leading coefficient of $p_k(x)$ is 1. The polynomials satisfy a three-term-recursion
\begin{align}\label{polrek}
x p_k(x) = p_{k+1}(x) + b_{k+1} p_k(x) + a_k p_{k-1}(x) \quad \mbox{ for } 1\leq k \leq n-1,
\end{align}
with $a_k >0$ and initial conditions $p_0(x) = 1,p_1(x) = x-b_1$. Indeed, Favard's Theorem (see \cite{szego1939}) says that any sequence of polynomials satisfying \eqref{polrek} with $a_k>0$ is orthogonal with respect 
to some measure $\mu$. Moreover, the measure $\mu$ is supported on $[0,\infty)$ if and only if there is a sequence of variables $z_k$ such that the decompositions
\begin{align}
a_k =& z_{2k-1} z_{2k} \label{zerlegungan},\\
b_k =& z_{2k-2} + z_{2k-1} \label{zerlegungbn}, 
\end{align}
hold with $z_k >0$ for $k \geq 1$ and $z_0=0$ (see \cite{chihara1978}). We will call the quantities $z_1,z_2,\dots $ recursion variables. If the support of the measure is contained in the compact interval $[0,1]$, the 
recursion variables form a chain sequence
\begin{align} \label{zerlegungzn}
 z_k =  (1-p_{k-1})p_k
\end{align}
with $p_k \in (0,1)$ for $k=1,\dots ,2n-1$ and $p_0 = 0$ (\cite{wall1948}). The $p_k$ are called canonical moments, where for a comprehensive overview on canonical moments we refer to the monograph of \cite{dettstud1997}. \\
Note that for any $k$, there is a one-to-one relation between the first $k$ moments of the measure $\mu$ and the first $k$ recursion coefficients or, if they exists, recursion variables or canonical moments. This can be seen by looking at the first orthogonal polynomials: fixing $2k-1$ moments is equivalent to fixing the first $k$ polynomials while the first $2k$ moments uniquely determine the first $k$ polynomials and additionally the norm of $p_k(x)$.\\
To be precise, the degree up to which the polynomials can be defined depends on the number of support points of $\mu$. If $a_1\dots a_k>0$, then orthogonal polynomials up to degree $k-1$ exist. In the cases we consider, all 
relevant polynomials can be defined. \\
 
\medskip

The limit eigenvalue distributions of the ensembles all have remarkably simple representations in the new variables: The recursion coefficients of the semicircle distribution satisfy $a_k=1, b_k=0$ for all $k$, the resulting polynomials are the Chebychev-polynomials of the second kind (rescaled on $[-2,2]$). For the Marchenko-Pastur law, the recursion variables $z_k$ are all equal to 1 and the canonical moments of the arcsine law are all given by $\tfrac{1}{2}$ (\cite{skibinsky1969}). 

\medskip

Now suppose that $\mu_n$ is the spectral measure of the tridiagonal matrix 
\begin{align} \label{tridiag}
T_n = \begin{pmatrix} 
  b_1 & \sqrt{a_1}    &         &         \\
                \sqrt{a_1} & b_2    & \ddots  &         \\
                    & \ddots & \ddots  & \sqrt{a_{n-1}} \\
                    &        & \sqrt{a_{n-1}} & b_n
\end{pmatrix}
\end{align}
with a positive subdiagonal, then the polynomials defined by the recursion coefficients $a_k,b_k$ are orthogonal with respect to $\mu_n$ (see \cite{deift2000}). In the language of the spectral theorem, $T_n$ represents the 
multiplication by the identity in the basis of orthogonal polynomials. \\

For the three classical ensembles, there exist tridiagonal reductions of the full matrix models such that for $\beta \in \{1,2,4\}$, the spectral measure of $T_n$ has the same distribution as the 
spectral measure of the full matrices. Moreover, the tridiagonal matrix yields the right distribution of eigenvalues and weights of $\mu_n$  for any value $\beta>0$. For the first two ensembles, the tridiagonal 
representation was obtained by \cite{dumede2002}. In the Gaussian case, 
\begin{align*}
a_k \sim \operatorname{Gamma}(\tfrac{\beta}{2}(n-k),\tfrac{2}{\beta n}), \qquad b_k \sim \mathcal{N}(0,\tfrac{2}{\beta n}) 
\end{align*}
and all recursion coefficients are independent. For the Laguerre ensemble with parameter $\gamma$, the spectral measure is supported on $[0,\infty)$ and the recursion variables satisfy
\begin{align*}
z_{2k-1} \sim \operatorname{Gamma}(\tfrac{\beta}{2}(n-k)+\gamma +1,\tfrac{2}{\beta n}), \qquad z_{2k} \sim \operatorname{Gamma}(\tfrac{\beta}{2}(n-k),\tfrac{2}{\beta n})
\end{align*}
and $z_1,\dots z_{2n-1}$ are independent. \cite{kilnen2004} constructed a reduction of the Jacobi ensemble. In this case, the canonical moments $p_1,\dots ,p_{2k-1}$ are independent and
\begin{align*}
p_k \sim \begin{cases} \operatorname{Beta}\left( \tfrac{2n-k}{4} \beta , \tfrac{2n-k-2}{4} \beta +\gamma+\delta+2 \right)\quad & k \mbox{ even,} \\
				\operatorname{Beta}\left( \tfrac{2n-k-1}{4} \beta +\gamma+1, \tfrac{2n-k-1}{4} \beta +\delta+1 \right)\quad & k \mbox{ odd.} \end{cases} 
\end{align*}
In any of the three cases, this provides us with a bijective mapping 
\begin{align*}
\psi_k: m^{(k)} \longmapsto r^{(k)}
\end{align*}
where $r^{(k)}$ is the vector of first $k$ independent variables (recursion coefficient, recursion variables or canonical moments, respectively). $\psi_k$ is continuous and if $\mu_n$ is a spectral measure, the new variables $r^{(k)}(\mu_n)$ satisfy $r^{(k)}(\mu_n) \to r^{(k)}(\sigma)$ almost surely. Moreover, we have (see \cite{detnag2012})
\begin{align*}
\frac{\partial \psi^{-1}_k}{\partial r^{(k)}} (r^{(k)}(\sigma)) = D_k
\end{align*}
with the matrix $D_k$ as in \eqref{clt}. 

\medskip

\section{Proofs}

\begin{prop}\label{mdps}
Let $\alpha,\beta,\gamma$ be positive constants, and $a_n$ be a sequence with $a_n\to \infty$ and $a_n/n\to 0$.
\begin{itemize}
\item[(i)] The sequence $(\sqrt{\tfrac{n}{a_n}}X_n)_n$, where $X_n\sim \mathcal{N}(0,\tfrac{\alpha}{n})$ satisfies a MDP with speed $a_n$ and good rate function $I_1(x) = \tfrac{1}{2\alpha}  x^2$.
\item[(ii)] The sequence $(\sqrt{\tfrac{n}{a_n}}(Y_n-1))_n$ with $Y_n \sim \operatorname{Gamma}(\alpha n+\beta,\tfrac{1}{\alpha n})$ satisfies a MDP with speed $a_n$ and good rate function $I_2(x) = \tfrac{\alpha}{2} x^2$.
\item[(iii)] The sequence $(\sqrt{\tfrac{n}{a_n}}(Z_n-\tfrac{1}{2}))_n$ with $Z_n \sim \operatorname{Beta}(\alpha n+\beta,\alpha n+\gamma)$ satisfies a MDP with speed $a_n$ and good rate function $I_3(x) = 4\alpha x^2$.
\end{itemize}
\end{prop}

\medskip

\proof
The first and second moderate deviations are well-known and follow from Theorem 3.7.1 in the book of \cite{demzei1998}. For the third sequence, we have the equality in distribution
\begin{align*}
Z_n = \frac{W_1}{W_1+W_2} =:\phi(W_1,W_2)
\end{align*}
with $W_1,W_2$ independent, $W_1 \sim \operatorname{Gamma}(\alpha n+\beta,\tfrac{1}{\alpha n}),\ W_2 \sim \operatorname{Gamma}(\alpha n+\gamma,\tfrac{1}{\alpha n})$. From $(ii)$ we know that the vector
\begin{align*}
\sqrt{n/a_n}(W_1-1,W_2-1)^T
\end{align*}
satisfies a MDP with speed $a_n$ and good rate function
\begin{align*}
I(x,y) = \tfrac{\alpha}{2}(x^2+y^2) .
\end{align*}
Then the $\delta$-method for moderate deviations (\cite{gaozhao2011}) implies a MDP for $\sqrt{\tfrac{n}{a_n}}(Z_n-\tfrac{1}{2}) = \sqrt{\tfrac{n}{a_n}}(\phi(W_1,W_2)-\phi(1,1))$ with good rate function
\begin{align*}
I_3(z) &= \inf \left\{ I(x,y)\ :\ \phi'(1,1)\cdot (x,y) = z\right\}\\
&= \inf \left\{ \tfrac{\alpha}{2}(x^2+y^2) :\ \tfrac{1}{4}x - \tfrac{1}{4}y = z\right\}\\
&= \tfrac{\alpha}{2} ((2z)^2+(-2z)^2) = 4{\alpha} x^2
\end{align*}
\QED

\medskip

\textbf{Proof of Theorem \ref{step1}:} \\
Let $\mu_n$ be the spectral measure of the Jacobi ensemble with parameters $\gamma, \delta$. Then by the tridiagonal reduction of the Jacobi ensemble in Section 3, the canonical moments $p_1,\dots ,p_{2n-1}$ are independent and 
\begin{align*}
p_k \sim \operatorname{Beta}(\tfrac{\beta}{2} n+c_k,\tfrac{\beta}{2} n+d_k) 
\end{align*}
with constants $c_k,d_k$ not depending on $n$. By Proposition \ref{mdps}, each normalised canonical moment $\sqrt{\tfrac{n}{a_n}}(p_k-\tfrac{1}{2})$ satisfies the MDP with speed $a_n$ and good rate function $I_3(x) = 2\beta x^2$. By the independence, the vector
\begin{align*}
\sqrt{\frac{n}{a_n}}(p_1-\tfrac{1}{2},\dots ,p_k-\tfrac{1}{2})^T
\end{align*}
satisfies the MDP with speed $a_n$ and rate $I_k(x) = 2\beta  ||x||^2$, $x\in \mathbb{R}^k$. Recall that $\psi^{-1}_k$ maps the first $k$ canonical moments to the first $k$ moments. The canonical moments of the arcsine law $\nu$ are all identical to $\tfrac{1}{2}$, so $\psi^{-1}_k(p_1,\dots ,p_k) = m^{(k)}(\mu_n)$ and $\psi^{-1}_k(\tfrac{1}{2},\dots ,\tfrac{1}{2}) = m^{(k)}(\nu)$. The $\delta$-method yields a MDP for the normalised moment vector $\sqrt{\tfrac{n}{a_n}}(m^{(k)}(\mu_n) - m^{(k)}(\nu))$ with good rate function
\begin{align*}
I_m(m) =& \inf \left\{ I_k(x)\ :\ \psi_k'(\tfrac{1}{2},\dots ,\tfrac{1}{2})\cdot x = m \right\} \\
=& 2 \beta \left| \left| \left( \psi_k'(\tfrac{1}{2},\dots ,\tfrac{1}{2})\right)^{-1} m \right| \right|^2 .
\end{align*}
The derivative of $\psi$ at $(\tfrac{1}{2},\dots ,\tfrac{1}{2})^T$ was first calculated by \cite{chakemstu1993} and is given by $2\sqrt{2} A$ with the $k\times k$ matrix $A$ given by \eqref{covjac}. This proves the theorem for the Jacobi ensemble. For the other two ensembles, the arguments are the same: first, Proposition \ref{mdps} implies a MDP for the independent auxiliary variables (the recursion coefficients or the recursion variables), then the $\delta$-method tranfers this MDP to the moments. The derivative appearing in the rate function is given by the matrices \eqref{covgau} and \eqref{covlag}, as calculated in the paper of \cite{detnag2012}. \QED

\medskip

\textbf{Proof of Lemma \ref{step2}:} \\
For the Jacobi ensemble, this is Lemma 2.1 in the paper of \cite{chakemstu1993}. In the Laguerre and Gaussian case, we use a combinatoric interpretation in terms of the generalised Catalan numbers
\begin{align*}
d_{i,j}= \binom{i+j}{i} - \binom{i+j}{i-1} .
\end{align*}
for $i\leq j$, see \cite{finu1976}. We remark that the Gaussian case could also be proven using the CLT-result of \cite{lytpast2009}, however since it follows directly from our arguments in the Laguerre case, we will follow this route. The generalised Catalan number $d_{i,j}$ counts the number of paths in the triangular lattice $\big( (i,j)\big)_{j\geq i}$ starting at $(i,j)$ and ending in $(0,0)$ and in each vertex are only allowed to move up or to the left (see Figure 1). For $i=j$ this is a usual Catalan path and we obtain the Catalan number $d_{i,i}$, which is the $i$-th moment of the Marchenko-Pastur distribution $\eta$ (\cite{nicspe2006}).

\psset{xunit=0.6cm,yunit=0.6cm,runit=1cm}
\begin{figure}[h] \label{path}
\centering
\begin{pspicture}(0,0)(8,7) 
\psline[linewidth=0.5pt](0,6)(7,6)
\psline[linewidth=0.5pt](1,5)(7,5)
\psline[linewidth=0.5pt](2,4)(7,4)
\psline[linewidth=0.5pt](3,3)(7,3)
\psline[linewidth=0.5pt](4,2)(7,2)
\psline[linewidth=0.5pt](5,1)(7,1)
\psline[linewidth=0.6pt, dotsep=2pt, linestyle=dotted](7,6)(7.5,6)
\psline[linewidth=0.6pt, dotsep=2pt, linestyle=dotted](7,5)(7.5,5)
\psline[linewidth=0.6pt, dotsep=2pt, linestyle=dotted](7,4)(7.5,4)
\psline[linewidth=0.6pt, dotsep=2pt, linestyle=dotted](7,3)(7.5,3)
\psline[linewidth=0.6pt, dotsep=2pt, linestyle=dotted](7,2)(7.5,2)
\psline[linewidth=0.6pt, dotsep=2pt, linestyle=dotted](7,1)(7.5,1)

\psline[linewidth=0.5pt](1,6)(1,5)
\psline[linewidth=0.5pt](2,6)(2,4)
\psline[linewidth=0.5pt](3,6)(3,3)
\psline[linewidth=0.5pt](4,6)(4,2)
\psline[linewidth=0.5pt](5,6)(5,1)
\psline[linewidth=0.5pt](6,6)(6,1)
\psline[linewidth=0.5pt](7,6)(7,1)
\psline[linewidth=0.6pt, dotsep=2pt, linestyle=dotted](6,1)(6,0.5)
\psline[linewidth=0.6pt, dotsep=2pt, linestyle=dotted](7,1)(7,0.5)

\psline[linewidth=2pt](7,1)(7,2)(6,2)(6,4)(4,4)(4,5)(1,5)(1,6)(0,6)

\psdots[dotscale=1.5](0,6)(7,1)
\rput[cb](0.4,6.5){\footnotesize $(0,0)$}
\rput[cb](7.8,1.4){\footnotesize $(5,7)$}

\end{pspicture}
\caption{A generalised Catalan path starting in $(5,7)$, thus contributing to $d_{5,7}$.}
\end{figure}
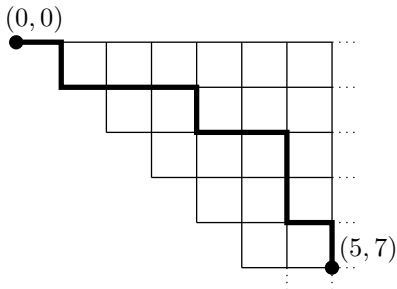

Recall that the matrix $D_k^{(L)}$ defining the covariance of the moments in the Laguerre case has entries
\begin{align*}
(D_k^{(L)})_{i,j} = \binom{2i}{i-j} - \binom{2i}{i-j-1} = d_{i-j,i+j} 
\end{align*}
for $j\ge i$ and 0 else. The asymptotic covariance of the $i$-th and the $j$-th moment is then given by
\begin{align*}
\left( D_k^{(L)}(D_k^{(L)})^T\right)_{i,j} =& \sum_{r=1}^{\min \{ i,j\} } d_{i-r,i+r} d_{j-r,j+r} \\
=& -d_{i,i} d_{j,j} + \sum_{r=0}^{\min \{ i,j\} } d_{i-r,i+r} d_{j-r,j+r}.
\end{align*}
Without loss of generality assume $i\leq j$. Consider a generalised Catalan path to $(0,0)$ in the triangular lattice starting in $(j-r,j+r)$. This path moves $j-r$ times up and $j+r$ times to the left. Replacing each move upwards by a move to the left and vice versa and reversing the order of movements, we obtain a unique path starting in $(i+j,i+j)$ and ending in $(i-r,i+r)$ (moving $j+r$ times up and $j-r$ times to the left). Therefore, 
\begin{align*}
 & \sum_{r=0}^i d_{i-r,i+r} d_{j-r,j+r} \\
 =& \sum_{r=0}^i \left( \# \left\{ \mbox{Paths from } (i-r,i+r) \mbox{ to } (0,0) \right\} \right) \left( \# \left\{ \mbox{Paths from } (i+j,i+j) \mbox{ to } (i-r,i+r) \right\} \right) .
\end{align*}
This last sum gives precisely the number of Catalan paths from $(i+j,i+j)$ to $(0,0)$ since each of those paths crosses the diagonal $\{(i-r,i+r)\}_{0\leq r\leq i}$ exactly once (see Figure 2). In conclusion,
\begin{align*}
\left( D_k^{(L)}(D_k^{(L)})^T\right)_{i,j} = d_{i+j,i+j} -d_{i,i} d_{j,j} = m_{i+j}(\eta) - m_i(\eta) m_j(\eta) .
\end{align*}

\psset{xunit=0.5cm,yunit=0.5cm,runit=0.5cm}
\begin{figure}[h] \label{path2}
\centering
\begin{pspicture}(0,0)(17,13)
\psline[linewidth=0.25pt](2,13)(15,13)
\psline[linewidth=0.25pt](3,13)(3,12)
\psline[linewidth=0.25pt](3,12)(15,12)
\psline[linewidth=0.25pt](4,13)(4,11)
\psline[linewidth=0.25pt](4,11)(15,11)
\psline[linewidth=0.25pt](5,13)(5,10)
\psline[linewidth=0.25pt](5,10)(15,10)
\psline[linewidth=0.25pt](6,13)(6,9)
\psline[linewidth=0.25pt](6,9)(15,9)
\psline[linewidth=0.25pt](7,13)(7,8)
\psline[linewidth=0.25pt](7,8)(15,8)
\psline[linewidth=0.25pt](8,13)(8,7)
\psline[linewidth=0.25pt](8,7)(15,7)
\psline[linewidth=0.25pt](9,13)(9,6)
\psline[linewidth=0.25pt](9,6)(15,6)
\psline[linewidth=0.25pt](10,13)(10,5)
\psline[linewidth=0.25pt](10,5)(15,5)
\psline[linewidth=0.25pt](11,13)(11,4)
\psline[linewidth=0.25pt](11,4)(15,4)
\psline[linewidth=0.25pt](12,13)(12,3)
\psline[linewidth=0.25pt](12,3)(15,3)
\psline[linewidth=0.25pt](13,13)(13,2)
\psline[linewidth=0.25pt](13,2)(15,2)
\psline[linewidth=0.25pt](14,13)(14,1)
\psline[linewidth=0.25pt](14,1)(15,1)
\psline[linewidth=0.25pt](15,13)(15,1)
\psline[linewidth=0.6pt, dotsep=2pt, linestyle=dotted](15,1)(15,0.5)

\psline[linewidth=0.6pt, linestyle=dashed](7,8)(12,13)

\psdots[dotscale=1](2,13)(10,11)(12,9)(14,1)(7,8)

\psline[linewidth=2.5pt, linestyle=dotted](2,13)(4,13)(4,12)(6,12)(6,11)(9,11)(9,9)(12,9)

\psline[linewidth=1.5pt](2,13)(7,13)(7,12)(10,12)(10,11)
\psline[linewidth=1.5pt](10,11)(10,8)(12,8)(12,5)(13,5)(13,3)(14,3)(14,1)

\rput[cb](13,0.4){\footnotesize $(i+j,i+j)$}
\rput[cb](6.5,7.4){\footnotesize $(i,i)$}
\rput[cb](2,12.4){\footnotesize $(0,0)$}

\end{pspicture}
\caption{A decomposition of a path from $(i+j,i+j)$ to $(0,0)$ along the (dashed) diagonal starting at $(i,i)$. The dotted line is a path starting in $(j-r,j+r)$, transformed into the path section from $(i+j,i+j)$ to $(i-r,i+r)$.}
\end{figure}
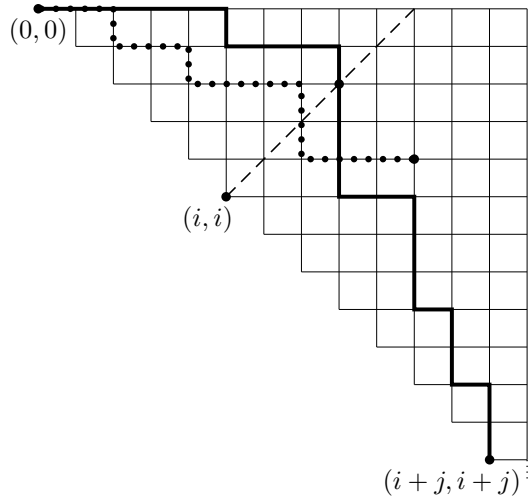

Now consider the remaining Gaussian ensemble, where the covariances are given by
\begin{align*}
\left( D_k^{(G)}(D_k^{(G)})^T\right)_{i,j} = \sum_{r=1}^{\min \{i,j\}} d_{\frac{i-r}{2},\frac{i+r}{2}} \mathbbm{1}_{\{ i+r \mbox{ even}\} } \cdot  d_{\frac{j-r}{2},\frac{j+r}{2}} \mathbbm{1}_{\{ j+r \mbox{ even}\} } .
\end{align*}
The even moments of the semicircle distribution $\rho$ are given by $m_{2k}(\rho) = d_{k,k}$ and the odd moments vanish.  
Suppose first that $i$ and $j$ are even, so we obtain from the calculations in the Laguerre case:
\begin{align*}
\left( D_k^{(G)}(D_k^{(G)})^T\right)_{i,j} =& - d_{\frac{i}{2},\frac{i}{2}} d_{\frac{j}{2},\frac{j}{2}} +  \sum_{r=0}^{\min \{i,j\}} d_{\frac{i-r}{2},\frac{i+r}{2}}  d_{\frac{j-r}{2},\frac{j+r}{2}} \mathbbm{1}_{\{ r \mbox{ even}\} } \\
=& - d_{\frac{i}{2},\frac{i}{2}} d_{\frac{j}{2},\frac{j}{2}} +  \sum_{r=0}^{\min \{i,j\}/2} d_{\frac{i}{2}-r,\frac{i}{2}+r}  d_{\frac{j}{2}-r,\frac{j}{2}+r}  \\
=& -d_{\frac{i}{2},\frac{i}{2}} d_{\frac{j}{2},\frac{j}{2}} + d_{\frac{i+j}{2},\frac{i+j}{2}} \\
=& m_{i+j}(\rho) - m_i(\rho) m_j(\rho) .
\end{align*}
If one of the indices $i,j$ is even and the other one is odd, all summands vanish and in this case also $m_i(\rho) m_j(\rho) - m_{i+j}(\rho)=0$. For the last case, assume $i,j$ are both odd, then
\begin{align*}
\left( D_k^{(G)}(D_k^{(G)})^T\right)_{i,j} =&   \sum_{r=0}^{\min \{i,j\}} d_{\frac{i-r}{2},\frac{i+r}{2}}  d_{\frac{j-r}{2},\frac{j+r}{2}} \mathbbm{1}_{\{ r \mbox{ odd}\} } \\
=&   \sum_{r=0}^{\min \{i-1,j-1\}/2} d_{\frac{i-1}{2}-r,\frac{i+1}{2}+r}  d_{\frac{j-1}{2}-r,\frac{j+1}{2}+r} \\
=&   \sum_{r=0}^{\min \{i^*,j^*\}} d_{i^*-r,i^*+1+r}  d_{j^*-r,j^*+1+r}
\end{align*}
with $i^*= \frac{i-1}{2}$ and $j^*=\frac{j-1}{2}$. Now assume that $i \leq j$ and apply the path transformation described above, where a generalised Catalan path starting in $(j^*-r,j^*+1+r)$ is uniquely identified with a path starting in $(i^*+j^*+1,i^*+j^*+1)$ and ending in $(i^*-r,i^*+1+r)$. Therefore,
\begin{align*}
& \sum_{r=0}^{\min \{i^*,j^*\}} d_{i^*-r,i^*+1+r}  d_{j^*-r,j^*+1+r} \\
=& \sum_{r=0}^{i^*} \left( \# \left\{ \mbox{Paths from } (i^*-r,i^*+1+r) \mbox{ to } (0,0) \right\} \right) \\ & \qquad \cdot  \left( \# \left\{ \mbox{Paths from } (i^*+j^*+1,i^*+j^*+1) \mbox{ to } (i^*-r,i^*+1+r) \right\} \right) \\
=& d_{i^*+j^*+1,i^*+j^*+1}
\end{align*}
The last identity follows since each path crosses the diagonal $\{ (i^*-r,i^*+1+r) \}_{r\geq0}$ exactly once. We conclude,
\begin{align*}
\left( D_k^{(G)}(D_k^{(G)})^T\right)_{i,j} = d_{\frac{i+j}{2},\frac{i+j}{2}} = m_{i+j}(\rho) = m_{i+j}(\rho) - m_i(\rho) m_j(\rho) .
\end{align*}
This proves the third case of the Lemma. \QED

\medskip

\textbf{Proof of Theorem \ref{haupt}:} \\
Let $\mathbb{M}_0$ be the set of all moment sequences $m= m(\mu)= (m_1(\mu),m_2(\mu),\dots)$ of signed measures $\mu \in \mathcal{M}_0$ with compact support, endowed with the topology of pointwise convergence and the Borel $\sigma$-algebra. 
Theorem \ref{step1} yields a MDP for the projection $m^{(k)}$ of $m$ onto the first $k$ coordinates. Denote by $D$ the infinite dimensional triangular matrix with upper left blocks given by $D_k$ as in Theorem \ref{step1}. By the Dawson-G\"artner Theorem (see \cite{demzei1998}) the random sequence
\begin{align*}
\sqrt{a_n^{-1}n}(m(\mu_n)-m(\sigma)) = m\big( \sqrt{a_n^{-1}n}(\mu_n - \sigma)\big) \in \mathbb{M}_0
\end{align*}
satisfies a moderate deviation principle with speed $a_n$ and good rate function
\begin{align*}
\mathcal{I}_m(m) = \sup_{k\in \mathbb{N}} \mathcal{I}_k(m^{(k)}) = \frac{\beta}{4} ||D^{-1} m||^2 .
\end{align*}
Note that since $D$ is a lower triangular matrix, $D^{-1} m$ is well-defined (with possible entries $\infty$). By the continuity of the mapping $m(\mu) \mapsto \mu$ from $\mathbb{M}_0$ to $\mathcal{M}_0$, the contraction principle yields a MDP for $\sqrt{a_n^{-1}n}(\mu_n - \sigma)$ with rate function 
\begin{align*}
\mathcal{I}(\mu) = \frac{\beta}{4} ||D^{-1} m(\mu)||^2 .
\end{align*}
To prove the integral representation of the rate function, we first apply an idea presented in the paper of \cite{chakemstu1993}. Define $X_i: x\mapsto x^i-m_i(\sigma)$ for $i\in \mathbb{N}$ as a random variable on $(\mathbb{R},\mathcal{B}(\mathbb{R}),\sigma)$, then
\begin{align*}
\mathbb{E}[X_i] = 0,\qquad \operatorname{Cov}(X_i,X_j) = m_{i+j}(\sigma) -m_i(\sigma)m_j(\sigma) .
\end{align*}
That is, Lemma \ref{step2} implies that the expectation and covariance of $X = (X_1,\dots ,X_k)^T$ coincides with the asymptotic expectation and covariance of $\sqrt{\frac{\beta n}{2}}  (m^{(k)}(\mu_n) - m^{(k)}(\sigma))$. But then the entries of $D_k^{-1}X$ are uncorrelated with variance 1, which means that $D_k^{-1}X$ has as entries the first $k$ polynomials $(p_1(x),\dots ,p_k(x))$ orthonormal with respect to $\sigma$. Therefore, the entries in the $i$-th row of $D_k^{-1}$ are given by the coefficients of the $i$-th orthonormal polynomial, except for the constant term (since $X_i$ is normalised). 
We can conclude that
\begin{align*}
D^{-1} m(\mu) & =  \left( \int p_1(x)-p_1(0) \mu(dx), \int p_2(x)-p_2(0) \mu(dx), \dots \right)^T \\ & = 
\left( \int p_1(x) \mu(dx), \int p_2(x) \mu(dx), \dots \right)^T
\end{align*}
where the last identity follows from $\mu(\mathbb{R})=0$. The rate function is then
\begin{align*}
\mathcal{I}(\mu) = \frac{\beta}{4} \sum_{k=1}^\infty \left( \int p_k(x) \mu(dx) \right)^2 
\end{align*}
and the assertion of the theorem follows from Lemma \ref{rate}. \QED



\medskip

\begin{lem} \label{rate}
Let $\sigma$ be a probability measure with infinite and compact support and let  $p_0,p_1,\dots$ denote the polynomials orthonormal with respect to $\sigma$. Then for any $\mu \in \mathcal{M}_0$, we have the equality of rate functions,
\begin{align*}
 \sum_{k=0}^\infty \left( \int p_k(x) \mu(dx) \right)^2 = \int \left(\tfrac{\partial \mu}{\partial \sigma} (x)\right)^2 \sigma (dx)
\end{align*}
and both sides equal $\infty$ if $\mu$ is not absolutely continuous with respect to $\sigma$.
\end{lem}

\medskip

\proof
First, note that if $d\mu = h\cdot d\sigma$ is absolutely continuous with respect to $\sigma $ with signed density $h$, Lemma \ref{rate} follows from Parseval's identity. We have to show that the left hand side equals $\infty$ as soon as $\mu$ is not absolutely continuous with respect to $\sigma$.

Suppose $\sigma(A)=0$ for a measurable set $A$, but $\mu(A)\neq 0$. We can choose polynomials $h_n$ of degree $n$ such that $\int x^k h_n(x)\sigma(dx)=\int x^k \mu(dx)$ for all $k\leq n$. Indeed, the infinite support of $\sigma$ implies that the Hankel-matrices $(m_{i+j}(\sigma))_{0\leq i,j\leq n}$ are positive definite (\cite{shotam1963}) and such polynomials exist. Then the measures $d\mu_n = h_n\cdot d\sigma $ have the first $n$ moments equal to the first $n$ moments of $\mu$ and therefore converge to $\mu$ in the moment topology. The support of $\mu_n$ is given by the support of $\sigma$ and is in particular uniformly bounded. Since $\int p_kd\mu$ depends only on the first $k$ moments of $\mu$, 
\begin{align*}
\sum_{k=0}^\infty \left( \int p_k(x) \mu(dx) \right)^2 = \lim_{n\to \infty} \sum_{k=0}^\infty \left( \int p_k(x) \mu_n(dx) \right)^2 = \lim_{n\to \infty} \int h_n(x)^2 \sigma (dx)
\end{align*}
Assume that there exists a subsequence of measures $\mu_{n_k}$ along which the integrals on the right hand side are uniformly bounded, i.e. $\sup_k  \int h_{n_k}(x)^2 \sigma (dx)<\infty$. Then the total variation 
\begin{align*}
\operatorname{TV} (\mu_{n_k}) = \sup_{f: |f|\leq 1} \int f d\mu_{n_k} \leq \sup_{f: |f|\leq 1} \left( \int (fh_{n_k})^2 d\sigma \right)^{1/2} \leq \left( \int h_{n_k}^2 d\sigma \right)^{1/2}
\end{align*}
is uniformly bounded as well. Let $f$ be a continuous function and for $\varepsilon>0$, let $p_f$ a polynomial approximating $f$ in the sup-norm on the compact support of $\sigma$ and $\mu$,
\begin{align*}
 \sup_{x\in \operatorname{supp}(\sigma)\cup  \operatorname{supp}(\mu)} |f(x)-p_f(x)| <\varepsilon ,
\end{align*}
then we can bound
\begin{align*}
\left| \int fd\mu_{n_k} - \int fd\mu \right| \leq \left| \int p_fd\mu_{n_k} - \int p_fd\mu \right|
+ \varepsilon (\operatorname{TV}(\mu_{n_k}) + \operatorname{TV}(\mu) ) .
\end{align*}
Since all moments converge, the right hand side can be made arbitrarily small as $n\to \infty$ and we get
$\mu_{n_k} \xrightarrow[n \rightarrow \infty ]{} \mu$ not only in the moment topology but also weakly. Moreover, the bounded total variation and Helly's selection theorem yield (possibly taking another subsequence) the weak convergence of the positive parts $\mu_{n_k}^+$ to the positive part $\mu^+$, when $\mu=\mu^+-\mu^-$ is the Hahn-Jordan decomposition of a measure $\mu$ into the nonnegative measures $\mu^+,\mu^-$. Without loss of generality, assume that the set $A$ is a subset of the support of $\mu^+$.
By the (outer) regularity of $\sigma$, for any $\varepsilon>0$ there is an open set $A_\varepsilon$ containing $A$ such that $\sigma(A_\varepsilon)<\varepsilon$. The Cauchy-Schwarz inequality implies
\begin{align*}
\int h_{n_k}^2 d\sigma \geq \sigma(A_\varepsilon)^{-1} \left( \int_{A_\varepsilon} |h_{n_k}| d\sigma  \right)^2
\geq  \varepsilon^{-1} \mu_{n_k}^+(A_\varepsilon)^2 
\end{align*}
Since $A_\varepsilon$ is open, we have $\liminf_{k\to \infty} \mu_{n_k}^+(A_\varepsilon)\geq \mu^+(A_\varepsilon)>0$ by the Portmanteau-Theorem and the weak convergence of $\mu_{n_k}^+$. Letting $\varepsilon \to 0$, this contradicts the boundedness of $\int h_{n_k}(x)^2 \sigma (dx)$ and we get $\lim_{n\to \infty} \int h_n^2 d\sigma =\infty$.
\QED

\medskip

\textbf{Acknowledgments.} The author would like to thank Fabrice Gamboa and Alain Rouault for several helpful discussions.

\bibliographystyle{apalike}
\bibliography{detnag}

\end{document}